\documentclass[11pt, a4paper,notitlepage]{article}

\usepackage{amssymb}
\usepackage{amsmath}
\usepackage{amsthm}
\usepackage{bm}

\newcommand{\Var}{\mathrm{Var}\,}
\newcommand{\eps}{\varepsilon}
\newcommand{\dd}{\mathrm{d}}

\newcommand{\RR}{\mathbb{R}}

\newcommand{\rr}{\mathbb{R}}

\newcommand{\bx}{\mathbf{x}}
\newcommand{\bu}{\mathbf{u}}
\newcommand{\bv}{\mathbf{v}}
\newcommand{\bw}{\mathbf{w}}

\newcommand{\bs}{\mathbf{s}}
\newcommand{\bk}{\mathbf{k}}
\newcommand{\bX}{\mathbf{X}}
\newcommand{\bW}{\mathbf{W}}
\newcommand{\bZ}{\mathbf{Z}}

\newcommand{\bbf}{\mathbf{f}}
\newcommand{\bsigma}{\bm{\sigma}}

\newcommand{\multint}{\int\cdots\int}

\newtheorem{thm}{Theorem}[section]
\newtheorem{lem}[thm]{Lemma}

\newtheorem{cor}[thm]{Corollary}
\theoremstyle{definition}
\newtheorem{rem}[thm]{Remark}
\newtheorem{ex}[thm]{Example}

\newtheorem{cond}[thm]{Condition}

\newenvironment{prf}{\begin{proof}[\bf Proof]}{\end{proof}}

\def\infint{\int_{-\infty}^\infty}
\def\ex{{\rm E\,}}
\def\var{{\rm Var\,}}

\newcommand{\cov}{{\rm Cov}\,}

\begin{document}

\title{Multivariate Nonparametric  Volatility Density Estimation}

\author{Bert van Es and Peter Spreij   \\[.5ex]
{\normalsize\sl Korteweg-de Vries Institute for Mathematics}\\
{\normalsize\sl Universiteit van Amsterdam}\\
{\normalsize\sl PO Box 94248}\\
{\normalsize\sl 1090GE Amsterdam}\\
{\normalsize\sl The Netherlands}}

\maketitle

\begin{abstract}
We consider a continuous-time stochastic volatility
model. The model contains a stationary volatility process, the
multivariate density of the finite dimensional distributions of
which  we aim to estimate. We assume that we observe the process
at discrete instants in time. The sampling times will   be
equidistant with vanishing distance.

A multivariate Fourier-type deconvolution kernel density estimator
based on the logarithm of the squared processes is proposed to
estimate the multivariate volatility density. An expansion of the
bias and a bound on the variance are derived.
\end{abstract}

\section{Introduction}

Let $S$ denote the log price process of some stock in a financial
 market. It is often assumed that $S$ can be modelled as the
solution of a stochastic differential equation or, more general,
as an It\^o diffusion process. So we assume that we can write
\begin{equation}\label{eq:s}
\dd S_t = b_t\, \dd t + \sigma_t\, \dd W_t, \ \ \ S_0=0,
\end{equation}
or, in integral form,
\begin{equation}\label{eq:si}
S_t = \int_0^t b_s\,\dd s + \int_0^t \sigma_s \,\dd W_s,
\end{equation}
where $W$ is a standard Brownian motion and the processes $b$ and
$\sigma$ are assumed to satisfy certain regularity conditions (see
Karatzas and Shreve (1991)) to have  the integrals in
(\ref{eq:si}) well-defined. In a financial context, the process
$\sigma$ is called the volatility process. One usually takes the
process $\sigma$ independent of the Brownian motion $W$.

In this paper we adopt this independence assumption and we model
$\sigma$ as a strictly stationary positive process satisfying a
mixing condition, for example an ergodic diffusion on
$(0,\infty)$. We will assume that all $p$-dimensional marginal
distributions of $\sigma$  have  invariant densities with respect
to the Lebesgue measure on $(0,\infty)^p$. This is typically the
case in virtually all stochastic volatility models that are
proposed in the literature, where the evolution of $\sigma$ is
modelled by a stochastic differential equation, mostly in terms of
$\sigma^2$, or $\log \sigma^2$ (cf.\ e.g.\ Wiggins (1987), Heston
(1993)).

As a motivation for nonparametric estimation procedures we
consider differential equations of the type
\[
\dd X_t = b(X_t)\, \dd t + a(X_t)\, \dd B_t,
\]
with $B$ equal to Brownian motion. Focussing on the invariant
univariate density, we recall that it is up to a multiplicative
constant equal to
\begin{equation}
\label{eq: piet}
x \mapsto \frac{1}{a^2(x)}\,{\exp \left(2\int_{x_0}^x
\frac{b(y)}{a^2(y)}dy\right)},
\end{equation}
where $x_0$ is an arbitrary element of the state space $(l,r)$, see
 e.g.\ Gihman and Skorohod (1972) or Skorokhod (1989). From
formula~(\ref{eq: piet}) one sees that the invariant distribution
of the volatility process (take $X$ for instance equal to
$\sigma^2$ or $\log \sigma^2$) may take  on many different forms,
as is the case for the various models that have been proposed in
the literature. Refraining from parametric assumptions on the functions $a$ and $b$,
nonparametric statistical procedures may be used to obtain information about the shape of the (one dimensional) invariant distribution.

A phenomenon that is often observed in practice,
is {\em volatility clustering}. This means that for different time
instants $t_1,\ldots,t_p$ that are close, the corresponding values
of $\sigma_{t_1},\ldots,\sigma_{t_p}$ are close again. This can
partly be explained by assumed continuity of the process $\sigma$,
but it might also result from specific areas where the
multivariate density of $(\sigma_{t_1},\ldots,\sigma_{t_p})$
assumes high values. For instance, it is conceivable that for
$p=2$, the density of $(\sigma_{t_1},\sigma_{t_2}$) has high
concentrations around points $(\ell,\ell)$ and $(h,h)$, with
$\ell<h$, a kind of bimodality on the diagonal of the joint distribution, with the
interpretation that clustering occurs around a low value $\ell$ or
around a high value $h$.

Here is an example where this happens. We consider a regime switching volatility process. Assume that for
$i=0,1$ we have two stationary processes $X^i$, each of them having
multivariate invariant distributions having densities. Call these
$f^i_{t_1,\ldots,t_p}(x_1,\ldots,x_p)$, whereas for $p=1$ we
simply write $f^i$. We assume these two processes to be
independent, and also independent of a two-state homogeneous
Markov chain $U$ with states $0,1$. Let $Q(t)$ be the matrix of
transition probabilities $q_{ij}(t)=P(X_t=i|X_0=j)$. Let $A$ be
the matrix of transition intensities and write
\[
A= \begin{pmatrix} -a_0 & a_1 \\
a_0 & -a_1
\end{pmatrix},
\]
with $a_0, a_1>0$. Then $\dot{Q}(t)=AQ(t)$, and
\[
Q(t)= \frac{1}{a_0+a_1}
\begin{pmatrix}
a_1+a_0e^{-(a_0+a_1)t} & a_1-a_1e^{-(a_0+a_1)t} \\
a_0-a_0e^{-(a_0+a_1)t} & a_0+a_1e^{-(a_0+a_1)t}
\end{pmatrix}.
\]
The stationary distribution of $U$ is given by
$\pi_i:=P(U_t=i)=\frac{a_{1-i}}{a_0+a_1}$ and we assume that $U_0$ has this distribution.
We finally define the process $\xi$ by
\[
\xi_t=U_tX^1_t+(1-U_t)X^0_t.
\]
Then $\xi$ is stationary too and it has a bivariate stationary
distribution with a density, related by $P(\xi_s\in\dd
x,\xi_t\in \dd y)=f_{s,t}(x,y)\,\dd x\,\dd y$.
Elementary calculations lead to the following expression for
$f_{s,t}$ for $0<s<t$.
\begin{align*}
f_{s,t}(x,y) & = q_{11}(t-s)\pi_1f^1_{s,t}(x,y) + q_{10}(t-s)\pi_0f^0(x)f^1(y) \\ & \mbox{} +q_{01}(t-s)\pi_1f^1(x)f^0(y) + q_{00}(t-s)\pi_0f^0_{s,t}(x,y).
\end{align*}
Suppose that the volatility process is defined by $\sigma_t=\exp(\xi_t)$ and that the $X^i$ are both Ornstein-Uhlenbeck processes given by
\[
\dd\, X^i_t = -a(X^i_t-\mu_i)\,\dd t + b\,\dd W^i_t,
\]
with $W^1$, $W^2$ independent Brownian motions, $\mu_1\neq \mu_2$ and $a>0$. Suppose that the $X^i$ start in their stationary $N(\mu_i,\frac{b^2}{2a})$ distributions. Then the centre of the distribution of $(X^i_s,X^i_t)$ is $(\mu_i,\mu_i)$, whereas the centre of the distribution of $(X^0_s,X^1_t)$ is $(\mu_0,\mu_1)$. Hence the density $f_{s,t}$ is a mixture of four hump shaped contours, each of them having a different centre of location. If $t-s$ is small, this effectively reduces to mixture of distributions with centres $(\mu_1,\mu_1)$ and $(\mu_2,\mu_2)$.

Nonparametric procedures are able to detect such
a property of a bivariate distribution, and are consequently by all means sensible
tools to get some partial insight in the behaviour of the
volatility.

In the present paper we propose a nonparametric  estimator for the
multivariate density of the volatility process. Using ideas from
deconvolution theory, we will propose a procedure for the
estimation of this  density at a number of fixed time instants.
Related work on estimating a univariate density has been done by
Van Es et al.~(2003), Comte and Genon-Catalot~(2006), Van Zanten
and Zareba~(2008), whereas a deconvolation approach has also been
adopted to estimate a regression function for a discrete time
stochastic volatility model  by Franke et al.~(2003), Comte~(2004)
and Comte et al.~(2008).

The observations of log-asset price $S$ process are assumed to
take place at the time instants $\Delta, 2\Delta, \ldots,
n\Delta$, where the time gap satisfies $\Delta = \Delta_n \to 0$
and $n\Delta_n \to \infty$ as $n \to \infty$.
This means that we base our estimator on so called {\em high frequency data}.

To asses the quality
of our procedure, we will study how the bias and variance of the
estimator behave under these assumptions. In Van Es et al.~(2003)
this problem has been studied for the marginal {\em univariate}
density of $\sigma$. The multivariate study of the present paper
largely builds on the approach of the cited paper, in particular
we will rely on a number of technical results that are contained
in it, but also we will borrow ideas from Van Es et al.~(2005), where a
multivariate problem for discrete time models has been studied.
Nevertheless, we will encounter a number of technical problems
that are not present in the univariate case, nor in the
multivariate case for discrete time models.

The remainder of the paper is organized as follows. In the next
Section~\ref{sec:construction}, we give the heuristic arguments
that motivate the definition of our estimator. In
Section~\ref{sec:results} the main results concerning the
asymptotic behaviour of the estimator are presented and discussed.
The proofs of the main theorems are given in
Section~\ref{sec:proofs}. They are based on a number of technical
lemmas, whose proofs are collected in
Section~\ref{sec:lemmaproofs}.

\section{Construction of the estimator}\label{sec:construction}

To motivate the construction of the estimator, we first consider~(\ref{eq:s})
without the drift term, so we
assume to have
\begin{equation*}
\dd S_t = \sigma_t \,\dd W_t, \ \ \ S_0=0.
\end{equation*}
It is  assumed that we observe the process $S$
at the discrete time instants $0$, $\Delta$,
$2\Delta, \ldots, n\Delta$.
For $i=1, 2, \ldots$ we work, as in
Genon-Catalot et al.\ (1998, 1999), with  the normalized
increments
\begin{equation}\label{eq:xdelta}
X^{\Delta}_i = \tfrac{1}{\sqrt{\Delta}}(S_{i\Delta} -
S_{(i-1)\Delta})= \tfrac{1}{\sqrt{\Delta}}
\int_{(i-1)\Delta}^{i\Delta}\sigma_t\,\dd W_t.
\end{equation}
For small $\Delta$, we have the  rough approximation
\begin{equation}\label{eq:xapprox}
X^{\Delta}_i  \approx
\sigma_{(i-1)\Delta}  \tfrac{1}{\sqrt{\Delta}}(W_{i\Delta} -
W_{(i-1)\Delta}) = \sigma_{(i-1)\Delta}Z^\Delta_i,
\end{equation}
where for $i=1, 2, \ldots$ we define
\[
Z^\Delta_i= \frac{1}{\sqrt{\Delta}}(W_{i\Delta} -
W_{(i-1)\Delta}).
\]
By the independence and stationarity of Brownian increments, the
sequence $Z^\Delta_1, Z^\Delta_2, \ldots$ is an i.i.d.\ sequence
of standard normal random variables. Moreover,  the sequence is
independent of the process $\sigma$ by assumption.

Let us first describe the  univariate
density estimator. Taking the logarithm of the square   of
 $X_i^\Delta$  we get
\begin{equation}
\log( (X^\Delta_i)^2)  \approx \log(
\sigma_{(i-1)\Delta}^2) + \log ( (Z_i^{\Delta})^2),
\end{equation}
where the terms in the sum are independent.
Assuming that
the approximation is sufficiently accurate we can use this approximate
convolution structure to
estimate the unknown
density $f$ of $\log(  \sigma_{i\Delta}^2)$ from the observed
$\log ((X^\Delta_i)^2)$.

Before we can define the estimator, we need some more notation.
Observe that the density of the `noise' $\log (Z_i^\Delta)^2$,
denoted by $k$, is given by
\begin{equation}\label{densityk}
k(x) = \frac{1}{\sqrt{2\pi}}\, e^{\tfrac{1}{2}x} e^{-
\tfrac{1}{2}e^{x}}.
\end{equation}
The characteristic function of the density $k$ is denoted by
$\phi_k$. We have $\phi_k(t)={1\over
\sqrt{\pi}}\,2^{it}\,\Gamma(\tfrac{1}{2}+it)$ and it's asymptotic
expansion $|\phi_k(t)|\  =\sqrt{2}\,e^{-\frac{1}{2}\pi
|t|}(1+O(\tfrac{1}{|t|}))$, for $|t|\to\infty$, see Van Es et
al.~(2003).

We will use a kernel function $w$, satisfying the following condition.
For examples of such kernels see Wand (1998).

\begin{cond}\label{cond:w} Let $w$ be a real symmetric function with real valued
symmetric characteristic function $\phi_w$ having support [-1,1].
Assume further
\begin{enumerate}
\item
$\infint |w(u)|du < \infty$
,
$\infint w(u)du=1$
,
$\infint u^2|w(u)|du<\infty$
,
\item
$\phi_w(1-t)=At^\rho + o(t^\rho),\quad\mbox{as}\ t\downarrow 0 $
for some $\rho >0$.
\end{enumerate}
\end{cond}

Following a well-known approach in statistical deconvolution
theory, we use a {\em deconvolution kernel density estimator}, see
e.g.\ Section 6.2.4 of Wand and Jones (1995). Having the
characteristic functions $\phi_k$ and $\phi_w$ at our disposal,
choosing a positive {\em bandwidth} $h$, we introduce the kernel
function
\begin{equation}\label{fourkernel}
v_h(x)={1\over 2\pi}\infint {{\phi_w(s)}\over\phi_k(s/h)}\
e^{-isx}ds
\end{equation}
and  the density estimator of the univariate density $f$ given by
\begin{equation}\label{eq:fnhc}
f_{nh}(x)={1\over nh}\sum_{j=1}^n v_h\left({{x-\log(
(X^\Delta_j)^2)}\over h}\right).
\end{equation}
One   easily verifies that the function $v_h$, and therefore also
the estimator $f_{nh}$, is real-valued. In Van Es et al.~(2003)
bias expansion and bounds on the variance of $f_{nh}(x)$ have been
obtained.

\medskip

In the present paper we will extend these results to a
multivariate setting, in which we will estimate the density
$f(\bx)=f_{t_1,\ldots,t_p}(\bx)$, with $\bx=(x_1,\ldots,x_p)$, of
a vector $\log (\sigma^2_{t_1},\ldots,\log \sigma^2_{t_p})$. Here
the $0<t_1<\ldots<t_p$ denote $p$ pre-specified time points. Below
we use boldface expressions for  (random) vectors. The expression
for the estimator of this density will be seen to be analogous to
the estimator in the univariate case, that has been analyzed in
Van Es et al.~(2003), and exhibits some similarity with the
estimator of a similar multivariate density in a discrete time framework as
treated in Van Es et al.~(2005).

What one ideally needs to estimate  $f(\bx)$ are observations of
$p$-dimensional random vectors that all have a density  equal to
$f$. This happens under the observation scheme that we have
introduced previously, if the $t_k$ are multiples of $\Delta$,
$t_k=i_k\Delta$ say. In that case, one should use
$(X^\Delta_{t_1+j\Delta},\ldots,X^\Delta_{t_p+j\Delta})$
for all the values of $j$ that are given by the observations. The
complicating factor is however, that the $t_k$ are not given as
multiples of $\Delta$, which on the other hand would lead to an
uninteresting estimation problem, if $\Delta\to 0$. Note also that
this kind of problem is not present, when one aims at estimating a
univariate marginal density of $\log \sigma^2_t$. All
$\log\sigma^2_t$, $t>0$ have the same marginal density.
%

We approach the problem as follows. Let us first introduce some
auxiliary notation. Write $(i_1^\Delta,\ldots,i_p^\Delta)$ for the
vector $([t_1/\Delta],\ldots,[t_p/\Delta])$ where $[ . ]$ denotes
the floor function.   We use $ \bX^\Delta_j $ to denote the
random vectors of lenght $p$
\[
\bX^\Delta_j=( X^\Delta_{j} ,\ldots,
X^\Delta_{i_p^\Delta-i_1^\Delta +j} ),\ j=1,\ldots,n-i_p^\Delta+i_1^\Delta.
\]
Hence it's $k$-th component is $X^\Delta_{i_k^\Delta-i_1^\Delta +j}$, $k=1,\ldots,p$.
Analogously, $\log ((\bX^\Delta_j)^2)$
denotes the vector
\[
\log (\bX^\Delta_j)^2=(\log (X^\Delta_{j})^2
,\ldots,\log (X^\Delta_{i_p^\Delta-i_1^\Delta +j})^2).
\]
Anywhere else in the sequel, we adhere to a similar notation.
Functions of a vector are assumed to be evaluated componentwise,
yielding again a vector.

Note that $\bX^\Delta_j$ is, by virtue of (\ref{eq:xapprox}),
approximately equal to the vector
\begin{equation}\label{eq:xtilde}
\tilde{\bX}_j:=(\sigma_{(j-1)\Delta}Z^\Delta_{j},\ldots,\sigma_{(i_p^\Delta -i_1^\Delta+j-1)\Delta)}Z^\Delta_{i_p^\Delta-i_1^\Delta +j})
\end{equation}
and that $(\log\sigma^2_{(j-1)\Delta},\ldots,\log\sigma^2_{(i_p^\Delta -i_1^\Delta+j-1)\Delta})$ has
density equal to $f_{i_1^\Delta\Delta,\ldots,i_p^\Delta\Delta}$
for every $j$, because of the assumed stationarity. Since
$\Delta\to 0$, one can expect that
$f_{i_1^\Delta\Delta,\ldots,i_p^\Delta\Delta}(\bx)\approx
f_{t_1,\ldots,t_p}(\bx)$. This motivates us to use the
observations $\bX_j^\Delta$, or rather the $\log
(\bX^\Delta_j)^2$, in the construction of a kernel estimator.

The kernel $\bw$ that we will use in the
multivariate case is just a product kernel,
$\bw(\bx)=\prod_{j=1}^pw(x_j)$. Likewise we take
$\bk(\bx)=\prod_{j=1}^p k(x_j)$ and the Fourier transforms
$\phi_\bw$ and $\phi_\bk$ factorize as well. Let $\bv_h$ be
defined by
\begin{equation}\label{eq:bvh}
\bv_h(\bx)={1\over (2\pi)^p}\int_{\mathbb{R}^p}
{{\phi_\bw(\bs)}\over\phi_\bk(\bs/h)}\ e^{-i\bs\cdot \bx}\,d\bs,
\end{equation}
where $\bs\in\mathbb{R}^p$ and $\cdot$ denotes inner product.
Notice that we also have the factorization
$\bv_h(\bx)=\prod_{j=1}^p v_h(x_j)$. \medskip\\
We finish this section by presenting the
multivariate density estimator $\bbf_{nh}(\bx)$ that we will use to estimate
$f(\bx)$. It is given by
\begin{equation}\label{eq:fnhp}
\bbf_{nh}(\bx)={1\over
(n-i_p^\Delta+i_1^\Delta)h^p}\sum_{j=1}^{n-i_p^\Delta+i_1^\Delta} \bv_h\Big({{\bx-\log
((\bX_j^{\Delta})^2)}\over h}\Big).
\end{equation}
Note that this estimator bears some similarity to, but also differs from the corresponding one for a
discrete time model  in Van Es et al.~(2005), where the
multivariate density of $(\sigma_{t+1},\ldots,\sigma_{t+p})$ at
consecutive time points is the object under study.

Under the assumption that the function $v_h(x)$ of (\ref{fourkernel}) integrates to one, an estimator of $f_{(t_1,\ldots,t_{p-1})}(x_1,\ldots,x_{p-1})$ is obtained by integrating out the variable $x_p$ in~(\ref{eq:fnhp}), which is of similar appearance. Further integration over the variables $x_2,\ldots,x_{p-1}$ reduces this estimator to the estimator of the univariate density given by~(\ref{eq:fnhc}) upon the substitution of $n$ by $n-i_p^\Delta+i_1^\Delta$.

\section{Results}\label{sec:results}

To derive the asymptotic behaviour of the estimator, we need a
mixing condition on the process $\sigma$.
For the sake of clarity,
we recall  the basic definitions. For a certain process
$X$ let ${\cal F}_a^b$ be the $\sigma$-algebra of events
generated by the random variables $X_t,\ a \le t \le b$. The
mixing coefficient $\alpha(t)$ is defined by
\begin{equation}\label{eq:alphamix}
\alpha(t)=\sup_{A\in {\cal F}_{-\infty}^0,\ B\in {\cal
F}_t^\infty} |P(A\cap B) -P(A)P(B)|.
\end{equation}
The process $X$ is called {\em  strongly mixing} if
$\alpha(t) \to 0$ as $t\to\infty$.
\medskip\\
As we mentioned in the introduction, it is common practice to
model the volatility process $V=\sigma^2$ as the stationary,
ergodic  solution of an SDE of the form
\[
d V_t = b(V_t)\,dt + a(V_t)\,d\,B_t.
\]
It is easily verified that for such processes it holds that $\ex
|V_t-V_0| = O(t^{1/2})$, provided that $b \in L_1(\mu)$ and $a \in
L_{2}(\mu)$, where $\mu$ is the invariant probability measure.
Indeed we have  $\ex |V_t-V_0| \leq \ex \int_0^t |b(V_s)|\,ds +
(\ex \int_0^t a^2(V_s)\, ds)^{1/2} = t||b||_{L_1(\mu)} + \sqrt{t}
||a||_{L_2(\mu)}$. In this setup, the process $V$ is strong
mixing, see for instance Corollary 2.1 of Genon-Catalot et
al.~(2000). Although we will not assume explicitly that $\sigma^2$
solves an SDE, the above observations motivate the the following
condition.
\begin{cond}\label{cond:sigma}
(i)  The process $\sigma$ is $L^1$-H\"older continuous of order
one half, $\ex|\sigma^2_t - \sigma_0^2| = O(t^{1/2})$ for $t \to
0$.
\\
(ii) The process $\sigma$ is strongly mixing with coefficient
$\alpha(t)$ satisfying,  for some $0<q<1$,
\begin{equation}\label{eq:alpha^q}
\int_0^\infty \alpha(t)^q\,dt<\infty.
\end{equation}
\end{cond}

\begin{rem}
Since the mixing coefficients $\alpha(t)$ are non-increasing in
$t$, condition~(\ref{eq:alpha^q}) is equivalent to the following.
For all  $t\in\RR$ there exists $C(q,t)$ such that for all
$\Delta>0$
\begin{equation}\label{eq:alpha^qt}
\sum_{k=0}^\infty \alpha(k\Delta + t)^q\leq \frac{C(q,t)}{\Delta},
\end{equation}
where $\alpha(t)$ is set equal to 1 for $t\leq 0$.
\end{rem}

\noindent
Our main theorems are multivariate versions of results in Van Es et
al.~(2003) which describe the asymptotic behaviour of the univariate
density estimator. Note that it also covers the case where there is
a drift $b_t$ present in equation (\ref{eq:s}). The condition on the
drift is boundedness of $\ex b_t^2$. This condition is typically
satisfied in realistic models for the log-returns of a stock, since
$b_t$ is the local rate of return and this will be mostly bounded
itself.

\begin{thm} \label{multcontasthm1}
Assume that $\ex b_t^2$ is bounded. Let the kernel function $w$
satisfy Condition~\ref{cond:w}. Let the density
$f_{t_1,\ldots,t_p}(\bx)$ of\/
$(\log\sigma^2_{t_1},\ldots,\log\sigma^2_{t_p})$ be continuous,
twice continuously differentiable with a bounded second derivative
and Lipschitz in $t_1,\ldots,t_p$, uniformly in $\bx$. Assume that
the first of Condition~\ref{cond:sigma} holds and that the
invariant density of $\sigma^2_t$ is bounded in a neighbourhood of
zero.
 Suppose that $\Delta=n^{-\delta}$ for  given $0<\delta <1$ and
choose $h=\gamma  \pi/\log n$, where $\gamma > 4p/\delta$. Then the
bias of the estimator (\ref{eq:fnhc}) satisfies
\begin{equation}\label{multcontasthm:1}
\ex \bbf_{nh}(\bx)= f_{t_1,\ldots,t_p}(\bx)+\tfrac{1}{2}h^2\int
\bu^\top \nabla^2f(\bx)\bu\, \bw(\bu)\,d\bu+o(h^2)+O(\Delta).
\end{equation}
\end{thm}
\begin{thm} \label{multcontasthm2}
Assume that $\ex b_t^2$ is bounded. Let the kernel function $w$
satisfy Condition~\ref{cond:w}. Assume that
Condition~\ref{cond:sigma} holds, that $\int |w(u)|^{2/(1-q)}\,du
<\infty$, where  $q$ is as in~(\ref{eq:alpha^q}), and that the
invariant density of $\sigma^2_t$ is bounded in a neighbourhood of
zero.
 Suppose that $\Delta=n^{-\delta}$ for  given $0<\delta <1$ and
choose $h=\gamma  \pi/\log n$, where $\gamma > 4p/\delta$. The
variance of the estimator satisfies
\begin{equation}\label{multcontasthm:2}
\var \bbf_{nh}(\bx) = O\Big({1\over
n}\,h^{2p\rho}e^{p\pi/h}\Big)+O\Big({1\over
{nh^{(1+q)p}\Delta}}\Big).
\end{equation}
\end{thm}

\begin{cor}
Under the assumptions of Theorems~\ref{multcontasthm1}
and~\ref{multcontasthm2} the bias satisfies  $\gamma^2\pi^2(\log n)^{-2}(1+o(1))$ and the order of the variance is $n^{-1+\delta}(\log
n)^{p(1+q)}$.
Hence the mean squared error of the estimator
$\bbf_{nh}(\bx)$ is of order $(\log n)^{-4}$.
\end{cor}
\begin{prf} The choices $\Delta=n^{-\delta}$, with $0<\delta<1$
and
$h=\gamma\pi/\log n$, with $\gamma>4p/\delta$ render a variance
that is of order $n^{-1+p/\gamma}(1/\log n)^{2p\rho}$ for the first
term of (\ref{multcontasthm:2}) and $n^{-1+\delta}(\log n)^{p(1+q)}$
for the second term. Since by assumption $\gamma>4p/\delta$ we
have $1/\gamma<\delta/4p<\delta$ so the second term dominates the
first term. The order of the variance is thus $n^{-1+\delta}(\log
n)^{p(1+q)}$. Of course, the order of the bias is logarithmic, hence
the bias dominates the variance and the mean squared error of
$f_{nh}(x)$ is also logarithmic.
\end{prf}
\noindent The proof of the theorems are deferred to the next
section. We conclude the present section by a number of comments
on the result.
\begin{rem}
The first order bound for the variance coincides with the order bound for the variance of the multivariate density estimator in discrete time models under the assumption that the volatility process and the error process are independent, see Theorem~3.2 in Van Es et al.~(2005). The second order bound is of the same nature as in the case of estimating a univariate density in continuous time models, see Theorem~3.1 in Van Es et al.~(2003), the difference being that in the multivariate case of the present paper one has $h^{p(1+q)}$ in the denominator instead of $h^{1+q}$.
\end{rem}

\begin{rem}
We observe some features that parallel
some findings for the univariate case. The expectation of the
deconvolution estimator is equal to the expectation of an ordinary
kernel density estimator, as becomes clear from the proof of
Lemma~\ref{lem:1}. It is well-known that the variance of
kernel-type deconvolution estimators heavily depends on the rate
of decay to zero of $|\phi_k(t)|$ as $|t|\to\infty$. The faster
the decay the larger the asymptotic variance.  This follows for
instance for i.i.d.\ observations from results in Fan (1991) and
for stationary observations from the work of Masry (1993). The
rate of decay of $|\phi_k(t)|$ for the density (\ref{densityk}) is
given by $ |\phi_k(t)|\ =\sqrt{2}\,e^{-\frac{1}{2}\pi
|t|}(1+O(\tfrac{1}{|t|})), $ see Lemma~5.3 in Van Es et
al.~(2003). This shows that $k$ is supersmooth, cf.~Fan (1991). By
the similarity of the tail of this characteristic function to the
tail of a Cauchy characteristic function  we can expect the same
order of the mean squared error as in Cauchy deconvolution
problems, where it decreases logarithmically in $n$, cf.~Fan
(1991) for results on i.i.d.\ observations. Note that this rate,
however slow, is faster than the one for normal deconvolution. Fan
(1991) also shows that we cannot expect anything better.
\end{rem}
\begin{rem}
The rate of convergence $(\log n)^{-4}$ for the mean squared error
as in Corollary 3.5 has also been found for other estimators.
Comte and Genon-Catalot~(2006) use (penalized) projection
estimators for $f$. These estimators are obtained by computing
certain projections on large but finite dimensional subspaces of
$L^2(\rr)$. Under similar assumptions as ours, they also find the
rate of convergence $(\log n)^{-4}$. By sharpening the assumed
smoothness properties of $f$, i.e. fast enough  exponential decay
of the characteristic function of $f$, so that $f$ itself is a
supersmooth density, they were able to obtain rates that are even
negative powers of $n$.

Van Zanten and Zareba~(2008) consider wavelet estimators of the
density of the accumulated squared volatility over intervals of
length $\Delta$ with $\Delta$ fixed for the model without drift
and with the same observations scheme. Under similar conditions,
they found this rate for the supremum of the mean integrated
squared error, the supremum taken over densities in some Sobolev
ball. For densities satisfying stronger smoothness conditions,
their estimators they obtained better rates, albeit still negative
powers of $\log n$. Both papers deal with estimating a univariate
density only.

Franke, H\"ardle and Krei\ss~(2003) consider a discrete time
model, where the evolution of $\log \sigma_t$ is decribed by a
nonlinear autoregression. By adopting a deconvolution approach
they estimate the unknown regression function and establish
tightness of the normalized estimators, where the normalization
again corresponds to the rate that we found.
\end{rem}

\begin{rem}
Better bounds on the asymptotic variance can be obtained under
stronger mixing conditions.
Consider for instance  {\em uniform mixing}. In this case the
mixing coefficient $\phi(t)$ is defined for $t >0$ as
\begin{equation}\label{eq:betamix}
\phi(t) =  \sup_{A\in {\cal F}_{-\infty}^0, B\in {\cal F}_t^{\infty}}
|P(A|B) -P(A)|
\end{equation}
and a process is called  {\em uniform
mixing} if $\phi(t)\to 0$ for $t\to \infty$. Obviously,  uniform
mixing implies strong mixing. As a matter of fact, one has the
relation
\[
\alpha(t) \leq \tfrac{1}{2}\phi(t).
\]
See Doukhan (1994) for this inequality and many other mixing
properties. If $\sigma$ is uniform mixing with coefficient
 $\phi$ satisfying
$\int_0^\infty \phi(t)^{1/2}dt<\infty$, then  the variance
bound is given by
\begin{equation}\label{contasthm:3}
\var f_{nh}(x) = O\Big({1\over
n}\,h^{2p\rho}e^{p\pi/h}\Big)+O\Big({1\over {nh^p\Delta}}\Big).
\end{equation}
The proof of this bound runs similarly to the strong-mixing bound.
The essential difference is that in equation~(\ref{deo}) we use
Theorem 17.2.3 of Ibragimov and Linnik (1971) with $\tau=0$  instead of Deo's (1973)
lemma, as in the proof of Theorem 2 in Masry (1983).
\end{rem}

\section{Proof of the Theorems}\label{sec:proofs}

We give the proof under the additional assumption that $b_t = 0$.
The general case is an easy consequence.

Let ${\cal F}_{\sigma}$ denote the sigma field generated by the
process $\sigma$. For $j=1,\dots, n-i_p^\Delta+i_1^\Delta$ we introduce,
along with the $\tilde{\bX}_j$ of~(\ref{eq:xtilde}), the following
vector notation
\begin{align*}
&\bsigma_j=( \sigma_{(j-1)\Delta} ,\ldots
\sigma_{(i_p^\Delta-i_1^\Delta +j-1)\Delta} )\\
&\bZ_j^\Delta =( Z_{j }^\Delta ,\ldots
Z_{i_p^\Delta-i_1^\Delta +j }^\Delta ),\\
\end{align*}
so that $\tilde{\bX}_j$ equals the Hadamard product $\bsigma_j\circ\bZ_j^\Delta$.
Note that since the $\sigma$ process is defined on the whole real
line the $\bsigma$ vectors are actually well defined for all $j$.\\
 Let
$\tilde \bbf_{nh}$ denote the estimator based on the approximating
random vectors $\tilde \bX_j$, i.e.
\begin{equation}
\tilde \bbf_{nh}(\bx)={1\over
(n-i_p^\Delta+i_1^\Delta)h^p}\sum_{j=1}^{n-i_p^\Delta+i_1^\Delta} \bv_h\Big({{\bx-\log
((\tilde \bX_j^{\Delta})^2)}\over h}\Big).
\end{equation}
The proof of (\ref{multcontasthm:1}) is partly based on the following two
lemmas, whose proofs are given in the next section. The first one
deals with the expectation  of $\tilde \bbf_{nh}$.

\begin{lem}\label{lem:1}
Let the density $f_{t_1,\ldots,t_p}(\bx)$ of\/
$(\log\sigma^2_{t_1},\ldots,\log\sigma^2_{t_p})$ be Lipschitz in
$t_1,\ldots,t_p$, uniformly in $\bx$. Then
\begin{equation}\label{eq:multexpf} \ex \tilde
\bbf_{nh}(\bx) = {1\over  h^p}\multint \bw\Big({{\bx-\bu}\over
h}\Big) f_{t_1,\ldots,t_p}(\bu)d\bu+O(\Delta)
\end{equation}
\end{lem}
Notice that, apart from the $O(\Delta)$ term, the
equality~(\ref{eq:multexpf}) is the same as for ordinary
multivariate kernel estimators, see for instance H\"ardle (1990) and
Scott (1992).
\medskip\\
The second lemma estimates the expected difference between $f_{nh}$ and $\tilde f_{nh}$.
The bound is in terms of the functions
\begin{equation}\label{eq: g0}
\gamma_0(h)={1\over 2\pi }\int_{-1}^1\Big|{{\phi_w(s)
}\over\phi_k(s/h)} \Big|ds
\end{equation}
and
\begin{equation}\label{eq: g1}
\gamma_1(h,x)=e^{\frac{1}{2}\pi/h} +
{1\over h}\exp\Big(\frac{\pi}{2}\frac{1+\pi/|x|}{h}\Big)\log\frac{1+\pi/|x|}{h}.
\end{equation}

\begin{lem}\label{expbound}
Assume Condition~\ref{cond:w} and that the first of
Condition~\ref{cond:sigma} holds and that the invariant density of
$\sigma^2_t$ is bounded in a neighbourhood of zero. For $h\to 0$
and $\eps$ small enough we have
\begin{align*}
\lefteqn{|\ex \bbf_{nh}(\bx) - \ex \tilde \bbf_{nh}(\bx)| =} \nonumber\\
&& O\left( {1\over h^{p+1}}\,
\gamma_0(h)^p{{\Delta}^{1/4}\over\eps} + {1\over
h^p}\,\gamma_0(h)^p\,{{\Delta^{1/2}}\over\eps^2} + {1\over
h^{p-1}}\gamma_0(h)^{p-1}\gamma_1(h,|\log 2\eps|/h){\eps\over|\log
2\eps|}\right).
\end{align*}
\end{lem}
\noindent {\bf Proof of Theorem~\ref{multcontasthm1}.} Statement
(\ref{multcontasthm:1}) follows by combining standard arguments of
kernel density estimation applied to
expression~(\ref{eq:multexpf}) in Lemma~\ref{lem:1} with
Lemma~\ref{expbound}. We will now show that the bound in
Lemma~\ref{expbound} is essentially a negative power of $n$,
whereas $h^2$ is of logarithmic order. Recall that we have assumed
$\delta > 4p/\gamma$. It follows that ${p}/{2\gamma} <
{\delta}/{4}-{p}/{2\gamma}$, so we can pick a $\beta \in
({p}/{2\gamma}, {\delta}/{4}-{p}/{2\gamma})$ and take
$\eps=n^{-\beta}$. By Lemmas \ref{gamma0bound} and \ref{thirdlem},
up to factors that are logarithmic in $n$, the order of $|\ex
\bbf_{nh}(\bx) - \ex \tilde \bbf_{nh}(\bx)|$ is then
\begin{equation}
n^{\frac{p}{2\gamma}-\frac{1}{4}\delta+\beta} +
n^{\frac{p}{2\gamma}+2\beta -\frac{\delta}{2}}
+n^{\frac{p}{2\gamma}-\beta},
\end{equation}
which is negligible to $h^2=\gamma^2 \pi^2/(\log n)^2$ for the
chosen values of the parameters.
\hfill$\square$\medskip\\
To prove the bound (\ref{multcontasthm:2}) we
 use the two lemmas below, which are proved in the next section. First  consider the
variance of $\tilde \bbf_{nh}(\bx)$.

\begin{lem}\label{multvar1}
Assume Condition~\ref{cond:w} and assume the second of
Condition~\ref{cond:sigma}. Assume also $\int |w(u)|^{2/(1-q)}\,du
<\infty$ for the same $q$ and $n\Delta\to\infty$. We
have, for $h \to 0$,
\begin{equation}\label{eq:twovar}
\var \tilde \bbf_{nh}(\bx) =O\Big({1 \over n } h^{2p\rho}
e^{p\pi/h}\Big) + O\Big({1\over {nh^{(1+q)p}\Delta}}\Big).
\end{equation}
\end{lem}
\noindent
The next lemma estimates $\var (f_{nh}(x)-\tilde f_{nh}(x))$.

\begin{lem}\label{multvarbound}
Assume that Condition~\ref{cond:w} and Condition~\ref{cond:sigma}
hold and let $\sigma^2_t$ have a bounded density in a
neighbourhood of zero. We have, for $h\to 0$ and $\eps > 0$ small
enough,
\begin{align}
& \var (\bbf_{nh}(\bx) -   \tilde \bbf_{nh}(\bx)) = \nonumber\\
& O\Big({1\over nh^{2p+2}} \gamma_0(h)^{2p}
{{\Delta}^{1/2}\over\eps^2}+{1\over nh^{2p-2}}
\gamma_0(h)^{2p-2}\,\gamma_1(h,|\log 2\eps|/h)^2{\eps\over|\log
2\eps|^2}\Big)\label{multvbound1a}\\
&+\frac{1}{nh^{2p}\Delta} O\Big({\Delta^{(1-q)/2}\over
{h^{2}\eps^{2}}}\, +\eps^{1-q }\Big).\label{multvbound2a}
\end{align}
\end{lem}
\begin{rem}
For $p=1$, the order bounds of Lemma~\ref{multvarbound} reduce to
those of Lemma~4.3 in Van Es et al.~(2003).
\end{rem}
\noindent {\bf Proof of Theorem~\ref{multcontasthm2}.} The bound
of (\ref{multcontasthm:2}) follows  as soon as we show that the
estimate in Lemma~\ref{multvarbound} is of lower order than the
one in Lemma~\ref{multvar1}. Up to terms that are logarithmic in
$n$, the bound in Lemma~\ref{multvar1} is of order $n^{\delta-1}$.
Choosing again $\eps=n^{-\beta}$, by Lemmas \ref{gamma0bound} and \ref{thirdlem},
one finds that, up to
logarithmic factors, the order of $\var( \bbf_{nh}(\bx) -  \tilde
\bbf_{nh}(\bx))$ is
\begin{equation}
n^{-1 +\frac{p}{\gamma}-\frac{\delta}{2}+2\beta} +n^{-1
+\frac{p}{\gamma}-\beta} + n^{-1+2\beta+\frac{1+q}{2}\delta} +
n^{-1+\delta-\beta(1-q)}.
\end{equation}
Recall our assumption $\delta\gamma > 4p$. If we pick $\beta$ less
than $\frac{1}{4}\,\delta(1-q)$, then all these terms are indeed of
lower order than $n^{\delta-1}$. \hfill$\square$

\section{Proof of  Lemmas \ref{lem:1}-\ref{multvarbound}}\label{sec:lemmaproofs}

We need expansions and order estimates for the functions $\phi_k$,
the kernel $v_h$ as defined in~(\ref{fourkernel}), $\gamma_0$ as
defined in (\ref{eq: g0}) and  the  function $\gamma_1$ as defined
in  (\ref{eq: g1}). These are collected in the next technical
lemmas, that are partially taken from Van Es et al.~(2003) and Van
Es et al.~(2005).


\begin{lem}\label{gamma0bound} Assume Condition~\ref{cond:w}. For $h\to 0$ we have
\begin{equation}
\gamma_0(h)=O\Big(h^{1+\rho}e^{\frac{1}{2}\pi /h}\Big).
\end{equation}
\end{lem}
\begin{prf}
See the proof of Lemma 5.3 in Van Es et al.~(2003).
\end{prf}

\bigskip

\begin{lem}\label{secondlem} Assume Condition~\ref{cond:w}. The functions $v_h$ and $\bv_h$ are bounded
and Lipschitz. More precisely, for all $x$  we have $|v_h(x)|\leq
\gamma_0(h)$ and for all $x$ and $u$ $|v_h(x+u)-v_h(x)|\leq
\gamma_0(h)\, |u|$. For all $p$ vectors $\bx$ we have
\begin{equation}\label{eq:bvbound}
\bv_h(\bx)|\leq \gamma_0(h)^p
\end{equation}
and for all $p$ vectors $\bx$ and $\bu$
\begin{equation}\label{multabscont}
|\bv_h(\bx+\bu)-\bv_h(\bx)|\leq \gamma_0(h)^p\, \sum_{j=1}^p|u_j|
\end{equation}
and for some $C>0$,
\begin{equation}\label{multabswcont}
|\bw(\bx+\bu)-\bw(\bx)|\leq  C \sum_{j=1}^p|u_j|.
\end{equation}
\end{lem}

\begin{prf}
The results for $|v_h(\cdot)|$  are known from Lemma 5.4 in Van Es
et al.~(2003). The bound~(\ref{eq:bvbound}) follows by the product
structure of $\bv_h$. Inequality~(\ref{multabscont}) follows by
induction and the same techniques can be used to prove
inequality~(\ref{multabswcont}).
\end{prf}

\begin{lem}\label{thirdlem} Assume Condition~\ref{cond:w}. For $x\to \infty$ we have the following
estimate on the behavior of $v_h$. For some positive constant $D$ it holds that
\begin{equation}\label{eq:vhx}
|v_h(x)|\leq  D {{\gamma_1(h,x)}\over{|x|}}\ \text{ as $
|x|\to\infty$},
\end{equation}
 and
\begin{equation}\label{gamma1bound}
\gamma_1(h,x)=O\Big({|\log h|\over h}\,e^{\frac{1}{2}\pi
(1+\pi/|x|)/h}\Big)\text{ as $h \to 0$.}
\end{equation}
Moreover, we have the following estimate on the behavior of
$\bv_h$. For some positive constant $D$ it holds that, if the
absolute value at least one of the components of $\bx$ tends to
infinity,
\begin{equation}\label{eq:vvhx}
|\bv_h(\bx)|\leq  D \gamma_0(h)^{p-1}{\gamma_1(h,x_1\vee\ldots\vee
x_p) \over{|x_1\vee\ldots\vee x_p|}} .
\end{equation}
\end{lem}
\begin{prf}
The estimates of (\ref{eq:vhx}) and (\ref{gamma1bound}) are taken
from Lemma 5,5 of Van Es et al.~(2003). To show (\ref{eq:vvhx}),
we argue as follows. Let $x^*=\max{x_1,\ldots,x_p}$. Without loss
of generality we may assume that $x^*=x_p$. Use the bound on $\gamma_0$ of
Lemma~\ref{secondlem} and the bound in~(\ref{eq:vhx}) to get
$|\bv_h(\bx)|=\prod_{i=1}^{p-1}v_h(x_i)v_h(x_p)\leq D
\gamma_0(h)^{p-1}\gamma_1(h,x_p)/x_p = D
\gamma_0(h)^{p-1}\gamma_1(h,x^*)/x^*$.
\end{prf}
\noindent
We are now ready with the proof of Lemma~\ref{lem:1}. Recall that ${\cal
F}_{\sigma}$ is the $\sigma$-algebra generated by the process
$\sigma$.

\bigskip
\noindent{\bf Proof of Lemma~\ref{lem:1}.} Write
\begin{eqnarray*}
\lefteqn{\ex (\tilde \bbf_{nh}(x)|{\cal F}_{\sigma})
=  {1\over
(n-i_p^\Delta+i_1^\Delta)h^p}\sum_{j=1}^{n-i_p^\Delta+i_1^\Delta} \ex
\Big(\bv_h\Big({{\bx-\log \bsigma^2_{ j- 1  } - \log
(\bZ^\Delta_j)^2 }\over h}\Big)|{\cal F}_{\sigma}\Big)}\\
&=& {1\over (n-i_p^\Delta+i_1^\Delta)h^p}\sum_{j=1}^{n-i_p^\Delta+i_1^\Delta}
{1\over (2\pi)^p}\multint{{\phi_\bw(\bs)}\over\phi_\bk(\bs/h)}\
\ex \Big(e^{-i\bs\cdot (\bx-\log \bsigma^2_{ j-1  } - \log
(\bZ^\Delta_j)^2)/ h}|{\cal F}_{\sigma}\Big)d\bs\\
& =& {1\over (n-i_p^\Delta+i_1^\Delta)h^p}\sum_{j=1}^{n-i_p^\Delta+i_1^\Delta}
{1\over (2\pi)^p}\multint {{\phi_\bw(\bs)}\over\phi_\bk(\bs/h)}\,
e^{-i\bs \cdot
\bx/h}e^{i\bs\cdot\log \bsigma^2_{ j-1  }/h}\,\phi_\bk(\bs/h)d\bs\\
&=& {1\over (n-i_p^\Delta+i_1^\Delta)h^p}\sum_{j=1}^{n-i_p^\Delta+i_1^\Delta}{1\over
(2\pi})^p\multint \phi_\bw(\bs)\ e^{-i\bs\cdot(\bx-\log
\bsigma^2_{ j-1 })/h}ds
\\
&=& {1\over (n-i_p^\Delta+i_1^\Delta)h^p}\sum_{j=1}^{n-i_p^\Delta+i_1^\Delta}
\bw\Big({{\bx-\log \bsigma^2_{ j-1  }}\over h}\Big).
\end{eqnarray*}
By taking the expectation we get, using
$|(i_j^\Delta-1)\Delta-t_j|\leq 2\Delta$, for $j=1,\ldots,p$, and the uniform Lipschitz continuity of $f$
\begin{align*}
\ex \tilde \bbf_{nh}&(x) = \ex \ex (\tilde \bbf_{nh}(x)|{\cal
F}_{\sigma})=\ex {1\over  h^p}\ex \bw\Big({{\bx-\log
\bsigma^2_{ 0  }}\over h}\Big)\\
& {1\over  h^p}\multint \bw\Big({{\bx-\bu}\over h}\Big)
f_{(i_1^\Delta-1)\Delta,\ldots,(i_p^\Delta-1)\Delta}(\bu)d\bu\\
&={1\over  h^p}\multint \bw\Big({{\bx-\bu}\over h}\Big)
f_{t_1,\ldots,t_p}(\bu)d\bu \\
& \quad+{1\over  h^p}\multint \bw\Big({{\bx-\bu}\over h}\Big)
(f_{(i_1^\Delta-1)\Delta,\ldots,(i_p^\Delta-1)\Delta}(\bu)-f_{t_1,\ldots,t_p}(\bu))d\bu\\
&={1\over  h^p}\multint \bw\Big({{\bx-\bu}\over h}\Big)
f_{t_1,\ldots,t_p}(\bu)d\bu+O(\Delta).
\end{align*}
\hfill$\square$
\medskip\\
For the proof of Lemma~\ref{expbound} we recall, see
Equations~(30) and~(31) in Van Es et al.~(2003), a few properties
of the process $\sigma$, valid under Condition~\ref{cond:sigma}.
There exists a constant $C> 0$ such that
\begin{equation}\label{sigma2}
\ex (X_1^\Delta -\sigma_{0}Z^\Delta_1)^2 \leq C \Delta^{1/2}
\text{ for $\Delta \to 0$},
\end{equation}
and
\begin{equation}\label{sigma3}
\ex \left|\tfrac{1}{\Delta}\int_0^\Delta \sigma^2_t\,dt -
\sigma_0^2 \right| \leq C\Delta^{1/2}\text{ for $\Delta \to 0$}.
\end{equation}
\medskip \\
{\bf Proof of Lemma~\ref{expbound}.}  We follow the lines of
thought as in the proof of Lemma 4.2 of Van Es et al.~(2003), now
applied in a multivariate setting. Let $\|\cdot\|$ denote the
Euclidean norm. Writing
\begin{equation}\label{multWjdef}
\bW_j=\bv_h\Big({{\bx-\log((\bX_j^\Delta)^2)}\over h}\Big)-
\bv_h\Big({{\bx-\log(\tilde \bX_j^2)}\over h}\Big),
\end{equation}
so that $\bbf_{nh}(\bx) -   \tilde \bbf_{nh}(\bx)={1\over
{(n-i_p^\Delta+i_1^\Delta)h^p}}\sum_{j=1}^{n-i_p^\Delta+i_1^\Delta} \bW_j$, and
defining the event $A$ as the event that all components of
$|\bX_1^\Delta|$ and $|\tilde \bX_1|$ are larger or equal to
$\eps$, we have
\begin{eqnarray}
|\ex \bbf_{nh}(\bx) - \ex \tilde \bbf_{nh}(\bx)| & \leq &
 {1\over h^p}\,\ex |\bW_1| \\
&=& {1\over h^p} \,\ex |\bW_1|
 I_A\label{multpartone}\\
&\quad& +{1\over h^p}\,\ex |\bW_1|
 I_{A^c} I_{[\|\bX_1^\Delta -\tilde \bX_1\|\geq\eps]}
\label{multparttwo}\\
&\quad& +{1\over h^p}\,\ex |\bW_1|
 I_{A^c} I_{[\|\bX_1^\Delta -\tilde \bX_1\|<\eps]}.
\label{multpartthree}
\end{eqnarray}
Recall that $|\log x-\log y|\leq |x-y|/\eps$ for $x,y\geq \eps$.
By Lemma~\ref{secondlem}, the bound~(\ref{sigma2}) and
stationarity, the term (\ref{multpartone}) can be bounded by
\begin{eqnarray*}
\lefteqn{{2\over h^{p+1}}\,\gamma_0(h)^p\sum_{j=1}^p\ex
|\log(|X_{i_j^\Delta}^\Delta|)-\log(|\tilde X_{i_j^\Delta}|)| I_A}\\
&\leq& {2p\over h^{p+1}}\,{1\over\eps} \gamma_0(h)^p\ex
|X_1^\Delta-\tilde X_1| \leq {2p\over h^{p+1}}\,
\gamma_0(h)^p\sqrt{C}{{\Delta}^{1/4}\over\eps}.
\end{eqnarray*}
This gives the first term in the order bound of
Lemma~\ref{expbound}.

The boundedness of the function $\bv_h$ as stated in
Lemma~\ref{secondlem} yields $|\bw_1|\leq 2\gamma_0(h)^p$. Using
also Chebychev's inequality and (\ref{sigma2}), we bound the term
(\ref{multparttwo}) by
\begin{align*}
{2\over h^p}\,\gamma_0(h)^pP(\|\bX_1^\Delta -\tilde
\bX_1\|\geq\eps)
& \leq {2\over h^p}\,\gamma_0(h)^pp P(|X_1^\Delta
-\tilde X_1|\geq {\eps\over \sqrt{p}}) \\
& \leq {2p^2\over
h^p}\,\gamma_0(h)^pC\,{\Delta^{1/2}\over\eps^2},
\end{align*}
which gives the second term in order bound of
Lemma~\ref{expbound}.

Consider the two arguments of the $\bv_h$ functions in $\bW_1$.
Since at least one of them  (and then the same for both arguments)
is in absolute value eventually larger than $|\log 2\eps |/h$, by
Lemma \ref{thirdlem} the term (\ref{multpartthree}) can be bounded
by \begin{align*}
2D{1\over h^p}&\gamma_0(h)^{p-1}\,\gamma_1(h,|\log
2\eps|/h){1\over(|\log 2\eps|/h)}\,pP(|\tilde X_1|\leq 2\eps)
\leq\\
&C_2 {1\over h^{p-1}}\gamma_0(h)^{p-1}\,\gamma_1(h,|\log
2\eps|/h){\eps\over|\log 2\eps|},
\end{align*}
for some constant $C_2$, where we used in the last inequality the
fact that the density of $\tilde{X}_1$ is bounded. This follows
from the assumption that
 $\sigma^2_0$ has a bounded density in a
neighbourhood of zero, as can easily be verified. \hfill$\square$

\bigskip

\noindent {\bf Proof of Lemma~\ref{multvar1}.} Consider the
decomposition
\begin{equation}\label{multvardecomp}
\var (\tilde \bbf_{nh}(\bx))= \var(\ex (\tilde \bbf_{nh}(\bx)|{\cal
F}_{\sigma})) + \ex(\var (\tilde \bbf_{nh}(\bx)|{\cal F}_{\sigma})).
\end{equation}
By the proof of Lemma \ref{lem:1} the conditional expectation $\ex
(\tilde \bbf_{nh}(\bx)|{\cal F}_{\sigma})$ is equal to a
multivariate kernel estimator of the density of $\log
\bsigma^2_{1}$. Adapting the proof of  Theorem 3 of Masry (1983) to
the multivariate situation, we can bound its variance by
$$
 {{20(1+o(1))}\over {nh^{(1+q)p}\Delta}}\ f(x)^{1-q}
\Big(\infint|w(u)|^{2/(1-q)}du \Big)^{1-q} \int_0^\infty
\alpha(\tau)^qd\tau,
$$
which is of the order $O(1/(nh^{(1+q)p}\Delta))$. This gives the
second order bound in~(\ref{eq:twovar}).

We turn to the expectation of the conditional variance.  Using
Lemma~\ref{secondlem}, we can bound the `diagonal terms' of the
conditional variance in~(\ref{multvardecomp})  by
$$
 {1\over
(n-i_p^\Delta+i_1^\Delta)h^{2p}}\, \ex \Big(\bv_h\Big({{\bx-\log
\tilde{\bX}^2_1}\over h}\Big)\Big)^2=O\Big({1\over
nh^{2p}}\gamma_0(h)^{2p}\Big),
$$
where we also used that $i^\Delta_p/n\to 0$.

Next we consider the `cross terms' of the conditional variance.
Since nonzero covariance can only occur if the vectors $\tilde
\bX_i$ and $\tilde \bX_j$ have common elements, we investigate a
`worst case'. For fixed $i$, there are at most $p-1$ among the
$\bx_j$ that have elements in common with $\bx_i$,
which yields
\begin{align*}
&{1\over (n-i_p^\Delta+i_1^\Delta)^2h^{2p}}\sum_{i\not=
j}\ex\cov\Big(\bv_h\Big({{\bx-\log \tilde{\bX}^2_i}\over
h}\Big),\bv_h\Big({{\bx-\log \tilde{\bX}^2_j}\over h}\Big)|{\cal
F}_\sigma\Big)\\
 &={2\over (n-i_p^\Delta+i_1^\Delta)^2h^{2p}}\sum_{i=1}^{n-i_p^\Delta+i_1^\Delta}\sum_{j=i+1}^{
i+i^\Delta_p- i^\Delta_1} \ex\cov\Big(\bv_h\Big({{\bx-\log
\tilde{\bX}^2_i}\over h}\Big),\bv_h\Big({{\bx-\log
\tilde{\bX}^2_j}\over h}\Big)|{\cal F}_\sigma\Big)\\
&\leq {2(p-1)\over (n-i_p^\Delta+i_1^\Delta)
h^{2p}}\gamma_0(h)^{2p}=O\Big({1\over n
h^{2p}}\gamma_0(h)^{2p}\Big),
\end{align*}
where in the last inequality we used that the expectation of the
conditional covariance is bounded in absolute value by $\ex
\Big(\bv_h\Big({{\bx-\log \tilde{\bX}^2_1}\over h}\Big)\Big)^2$,
due to stationarity. The first order bound in~(\ref{eq:twovar})
follows by an application of Lemma~\ref{gamma0bound}.
\hfill$\square$\medskip\\
{\bf Proof of Lemma~\ref{multvarbound}.} We will use arguments
similar to those in the proof of Lemma~\ref{expbound}. With
$\bW_j$ as in (\ref{multWjdef}) we have, using the ordinary
variance decomposition and stationarity of the $\bW_j$,
\begin{align}
\lefteqn{ \var( \bbf_{nh}(\bx) -   \tilde \bbf_{nh}(\bx))
 =} \nonumber \\ && \mbox{} {1\over (n-i_p^\Delta+i_1^\Delta) h^{2p}}
\var \bW_1 + {1\over (n-i_p^\Delta+i_1^\Delta)^2 h^{2p}}\sum_{i\not= j} \cov( \bW_i ,
\bW_j ).\label{eq:varcov}
\end{align}
Let us first derive a bound on $ \Var \bW_1$. As in the proof of
Lemma~\ref{expbound} we use  $A$, the event that all components of
$|\bX_1^\Delta|$ and $|\tilde \bX_1|$ are larger than or equal to
$\epsilon$. We have $ \Var \bW_1 \le  \ex \bW_1^2$, which can be
split up as the three terms sum
\begin{align}
\ex \bW_1^2 & = \ex \bW_1^2
 I_A\label{multvarpartone}\\
& \quad +  \ex \bW_1^2
 I_{A^c} I_{[\|\bX_1^\Delta -\tilde \bX_1\|\geq\eps]}
\label{multvarparttwo}\\
&\quad +  \ex \bW_1^2
 I_{A^c} I_{[\|\bX_1^\Delta -\tilde \bX_1\|<\eps]}.
\label{multvarpartthree}
\end{align}
By stationarity,  the Lipschitz property of $\bv_h$ in
Lemma~\ref{secondlem} and (\ref{sigma2}) the term
(\ref{multvarpartone}) can be bounded by
\begin{eqnarray}
\lefteqn{{4\over h^2}\,\gamma_0(h)^{2p}\ex\Big(\sum_{j=1}^p
|\log |X_{i_j^\Delta}^\Delta|-\log |\tilde X_{i_j^\Delta}||\Big)^2 I_A}\nonumber \\
&\leq&{4p\over h^2}\,\gamma_0(h)^{2p}\ex \sum_{j=1}^p
(\log|X_{i_j^\Delta}^\Delta|-\log|\tilde X_{i_j^\Delta}|)^2
I_A \nonumber \\
&\leq& {4p^2\over h^2}\,{1\over\eps^2} \gamma_0(h)^{2p}\ex
(|X_1^\Delta|-|\tilde X_1|)^2 \nonumber \\
&\leq& {4p^2\over h^2}\,{1\over\eps^2} \gamma_0(h)^{2p}\ex
(X_1^\Delta-\tilde X_1)^2 \leq {4p^2 \over h^2} \gamma_0(h)^{2p} C
{{\Delta}^{1/2}\over\eps^2}.\label{eq:varw1a}
\end{eqnarray}
We turn to the term (\ref{multvarparttwo}). By the  bound on
$\bv_h$ of Lemma~\ref{secondlem} and  by (\ref{sigma2}) again, it
can be bounded by
\begin{align*}
4\gamma_0(h)^{2p}P(\|\bX_1^\Delta -\tilde \bX_1\|\geq\eps)
& \le 4\gamma_0(h)^{2p}p P(|X_1^\Delta -\tilde X_1|\geq {\eps\over
\sqrt{p}}) \\
& \leq 4p^2\gamma_0(h)^{2p}C\,{\Delta^{1/2}\over\eps^2}.
\end{align*}
Due to absence of a factor $h^2$ in the denominator, this bound is
of smaller order than the one for~(\ref{multvarpartone}) and will
therefore be neglected.

Next we consider (\ref{multvarpartthree}). Recall form the proof
of Lemma~\ref{expbound} that $P(|\tilde X_1|\leq 2\eps)=O(\eps)$.
Since  at least one (the same) coordinate of the absolute value of
both arguments of $\bv_h$   is
 eventually larger than $|\log 2\eps |/h$,
by Lemma \ref{thirdlem} the term (\ref{multvarpartthree}) can be
bounded by
\begin{align}
4D^2 &\gamma_0(h)^{2p-2}\,\gamma_1(h,|\log 2\eps|/h)^2{1\over(|\log
2\eps|^2/h^2)}\,pP(|\tilde X_1|\leq
2\eps) \leq \nonumber \\
& C_2   h^2 \gamma_0(h)^{2p-2}\,\gamma_1(h,|\log
2\eps|/h)^2{\eps\over|\log 2\eps|^2}, \label{eq:varw1b}
\end{align}
for some constant $C_2$.

Wrapping up the order bounds~(\ref{eq:varw1a})
and~(\ref{eq:varw1b}) for $\ex \bW_1^2$, we get
\begin{equation}\label{multexpz}
\ex \bW_1^2 = O\Big( {1\over h^2} \gamma_0(h)^{2p}
{{\Delta}^{1/2}\over\eps^2}+h^2
\gamma_0(h)^{2p-2}\,\gamma_1(h,|\log 2\eps|/h)^2{\eps\over|\log
2\eps|^2}\Big),
\end{equation}
which, substituted in (\ref{eq:varcov}), gives the order bounds
of~(\ref{multvbound1a}).

We now consider the covariance terms in~(\ref{eq:varcov}), that
will be seen to have the order bounds of~(\ref{multvbound2a}). We
have the decomposition
\begin{equation}\label{eq:ecove}
\cov( \bW_i , \bW_j )=\ex\cov(\bW_i,\bW_j|{\cal F}_\sigma) +
 \cov(\ex(\bW_i|{\cal F}_\sigma),\ex(\bW_j|{\cal F}_\sigma)).
\end{equation}
The last term in~(\ref{eq:varcov}) then becomes
\begin{align}
\lefteqn{{2\over (n-i_p^\Delta+i_1^\Delta)^2
h^{2p}}\sum_{i=1}^{n-i_p^\Delta+i_1^\Delta}\sum_{j\not=i}^{ i+i^\Delta_p-
i^\Delta_1} \ex\cov(\bW_i,\bW_j|{\cal F}_\sigma)\label{term2}} \\
&\quad +{1\over (n-i_p^\Delta+i_1^\Delta)^2 h^{2p}}\sum_{i\not= j}
\cov(\ex(\bW_i|{\cal F}_\sigma),\ex(\bW_j|{\cal
F}_\sigma)).\label{term3}
\end{align}
In a first step we consider the expectation of the conditional
covariances in~(\ref{term2}). Arguing as in the proof of
Lemma~\ref{multvar1}, we can bound it by
\[
{(p-1)\over (n-i_p^\Delta+i_1^\Delta) h^{2p}}\var \bW_1,
\]
which is $p-1$ times the first term on the right hand side of
Equation~(\ref{eq:varcov}).  Hence its contribution can be
absorbed in the already obtained bounds of~(\ref{multvbound1a}).

%

\noindent Next we concentrate on the sum of covariances in
(\ref{term3}). Define
\begin{equation}
\bar{\sigma_i} = {1\over\Delta}\int_{(i-1)\Delta}^{i\Delta}
\sigma_t^2dt
\end{equation}
and the vector $\bar\bsigma_j$ by
$\bar\bsigma_j=(\bar\bsigma_{i_1^\Delta+j-1},\ldots,\bar\bsigma_{i_p^\Delta+j-1})$.
Note that given ${\cal F}_\sigma$, $\bX_i^\Delta$ is a
multivariate normal vector with independent components with
variances equal to the components of $\bar\bsigma_i$ and that
$\tilde{\bX_i}$ is a multivariate normal vector with independent
components with variances equal to the components of
$\bsigma_{i-1}^2$. As in the proof of Lemma~\ref{lem:1} it follows
that
$$
\ex( \bW_i|{\cal F}_\sigma)=\bw\Big({{\bx-\log \bar\bsigma_i}\over
h}\Big) -\bw\Big({{\bx-\log \bsigma_{(i-1) }^2}\over h}\Big).
$$

We follow the line of arguments in the proof of Theorem 3 in Masry
(1983). The stationarity of $\bW_j$ implies that also the
conditional expectations $\tilde{\bW}_j:=\ex(\bW_j|{\cal F}_\sigma)$
are stationary. Hence we have
\[
\sum_{i\not= j}\cov(\tilde{\bW}_i,\tilde{\bW}_j)= 2
\sum_{k=1}^{n-1}(n-k)\cov(\tilde{\bW}_0,\tilde{\bW}_k).
\]
Now note that the process $\tilde{\bW}_j$ is strongly mixing with
a mixing coefficient
$\tilde\alpha (k)\leq \alpha( (k-2) \Delta + t_1-t_p),
k=1,2,\dots$ if $k \Delta > t_p-t_1+2\Delta$ and
$\tilde\alpha(k)=1$ else. By a lemma of Deo (1973) for strongly
mixing processes it follows that for all $\tau>0$
\begin{equation}\label{deo}
|\cov(\tilde{\bW}_0,\tilde{\bW}_k)| \le 10\alpha((k-2) \Delta +
t_1-t_p)^{\tau/(2+\tau)}
\Big(\ex|\tilde{\bW}_1|^{2+\tau}\Big)^{{2/(2+\tau)}}.
\end{equation}
By the equivalent Condition~(\ref{eq:alpha^qt}) on the mixing
coefficients $\alpha(t)$ (applied with $\tau=2q/(1-q)$, a choice
for $\tau$ that we will make later on as well), we get
for~(\ref{term3})
\begin{eqnarray*}
\lefteqn{\Big|{1\over (n-i_p^\Delta+i_1^\Delta)^2 h^{2p}}\sum_{i\not= j}
\cov(\tilde{\bW}_i,
\tilde{\bW}_j)\Big|}\\
&\leq& {10\over
(n-i_p^\Delta+i_1^\Delta)h^{2p}}\,\Big(\ex|\tilde{\bW}_1|^{2+\tau}\Big)^{{2/(2+\tau)}}
 \sum_{k=1}^{n-i_p^\Delta }(1-{k\over
n-i_p^\Delta+i_1^\Delta})\alpha(k \Delta +
t_1-t_p)^{\tau/(2+\tau)}\\
&\leq& {10\over (n-i_p^\Delta+i_1^\Delta)h^{2p}}
\frac{C(\frac{\tau}{2+\tau},t_1-t_p-2\Delta)}{\Delta}
\Big(\ex|\tilde{\bW}_1|^{2+\tau}\Big)^{{2/(2+\tau)}}.
\end{eqnarray*}
Next we derive a bound on  $\ex|\tilde{\bW}_1|^{2+\tau}$. Fix
$\kappa\in (0,1]$ and define the event  $B$ as the event that all
components of $|\bar\bsigma_1|$ and $|\bsigma_0^2 |$ are larger or
equal to $\epsilon$. We have
\begin{eqnarray}
\lefteqn{\ex|\tilde{\bW}_1|^{2+\tau} = \ex
\Big|\bw\Big({{\bx-\log(\bar\bsigma_1)}\over h}\Big)-
\bw\Big({{\bx-\log(\bsigma_{0}^2)}\over h}\Big)\Big|^{2+\tau}I_B}
\label{multvarwpartone}\\
&+& \ex \Big|\bw\Big({{\bx-\log(\bar\bsigma_1)}\over h}\Big)-
\bw\Big({{\bx-\log(\bsigma_{0}^2)}\over h}\Big)\Big|^{2+\tau}
I_{B^c} I_{[\|\bar\bsigma_1^\kappa
-\sigma_{0}^{2\kappa}\|\geq\eps]}
\label{multvarwparttwo}\\
&+& \ex \Big|\bw\Big({{\bx-\log(\bar\bsigma_1)}\over h}\Big)-
\bw\Big({{\bx-\log(\bsigma_{0}^2)}\over h}\Big)\Big|^{2+\tau}
I_{B^c} I_{[\|\bar\bsigma_1^\kappa
-\bsigma_{0}^{2\kappa}\|<\eps]}. \label{multvarwpartthree}
\end{eqnarray}
By Lemma~\ref{secondlem} the term (\ref{multvarwpartone}) can be
bounded by a constant times
\begin{eqnarray}
\lefteqn{\frac{1}{h^{2+\tau}}\ex\Big(\sum_{j=1}^p|\log(\bar\sigma_{i_j^\Delta})-\log(
\sigma_{i_j^\Delta-1}^2)|\Big)^{2+\tau}I_B} \nonumber\\
&\leq&\frac{p^{1+\tau}}{h^{2+\tau}}\ex
\sum_{j=1}^p|\log(\bar\sigma_{i_j^\Delta})-\log(
\sigma_{i_j^\Delta-1}^2)|^{2+\tau}I_B\nonumber\\
&\leq&\frac{p^{2+\tau}}{h^{2+\tau}}\ex
 |\log(\bar\sigma_{i_1^\Delta})-\log(
\sigma_{i_1^\Delta-1}^2)|^{2+\tau}I_B\nonumber\\
&\leq&\frac{p^{2+\tau}}{(\kappa\epsilon h)^{2+\tau}}\ex
 | \bar\sigma_{1}^\kappa-
\sigma_{0}^{2\kappa} |^{2+\tau}.\label{eq:bound1}
\end{eqnarray}
The term (\ref{multvarwparttwo}) can be bounded by
\begin{equation*}
 pP(|\bar\sigma_1^\kappa -\sigma_{0}^{2\kappa}|\geq\frac{\eps}{\sqrt{p}})
\leq {p^{2+\tau/2}\over\eps^{2+\tau}}\,\ex
|\bar\sigma_1^\kappa-\sigma_{0}^{2\kappa}|^{2+\tau}.
\end{equation*}
Since this is for $h\to 0$ of smaller order
than~(\ref{eq:bound1}), it will be neglected in the sequel.

Finally we analyze  the term (\ref{multvarwpartthree}). On the
complement of $B$ there is at least one component of either
$|\bar\bsigma_1|$ or $|\bsigma_0^2 |$ that is smaller  or equal to
$\epsilon$. Together with $\|\bar\bsigma_1^\kappa
-\bsigma_{0}^{2\kappa}\|<\eps$ this implies that there is at least
one pair of corresponding components of the vectors that are both
smaller than $\eps(1+\eps^{1-\kappa})^{1/\kappa}$. Using the
stationarity, we bound the term (\ref{multvarwpartthree}) by
$$
pP(\bar\sigma_1\leq
\eps(1+\eps^{1-\kappa})^{1/\kappa}\ \mbox{and}\ \sigma_{0}^2\leq
\eps(1+\eps^{1-\kappa})^{1/\kappa}),
$$
which is bounded by
\begin{equation}\label{eq:bound2}
pP(\sigma_{0}^2\leq 2\eps)=O(\eps),
\end{equation}
since
$\sigma^2_0$ was assumed to have a bounded density in a
neighbourhood of zero.
\\
Combining ~(\ref{eq:bound1}) and~(\ref{eq:bound2}) with
$\tau=2q/(1-q)$ and $\kappa=\frac{1}{2+\tau}=\frac{1-q}{2}$,
 we have  with an application of the basic inequality
$|u^\kappa-v^\kappa|\leq |u-v|^\kappa$ for $u, v \geq 0$ and
$\kappa\in(0,1]$ in the second equality below and~(\ref{sigma3})
in the fourth equality for the term~(\ref{term3})
\begin{eqnarray*}
\lefteqn{\Big|{1\over (n-i_p^\Delta+i_1^\Delta)^2 h^{2p}}\sum_{i\not= j}
\cov(\tilde{\bW}_i,
\tilde{\bW}_j)\Big|}\nonumber\\
&=& \frac{1}{(n-i_p^\Delta+i_1^\Delta)h^{2p}\Delta}O\Big({1\over
h^{2+\tau}}\,{1\over\eps^{2+\tau}}\, \ex
|\bar\sigma_1^\kappa-\sigma_{0}^{2\kappa}|^{2+\tau}
+\eps\Big)^{2/(2+\tau)}\nonumber\\
&=&  \frac{1}{(n-i_p^\Delta+i_1^\Delta)h^{2p}\Delta}O\Big({1\over
h^{2+\tau}}\,{1\over\eps^{2+\tau}}\, \ex
|\bar\sigma_1-\sigma_{0}^{2}|^{\kappa(2+\tau)}
+\eps\Big)^{2/(2+\tau)}  \nonumber\\
&=& \frac{1}{(n-i_p^\Delta+i_1^\Delta)h^{2p}\Delta} O\Big({(\ex
|\bar\sigma_1-\sigma_{0}^{2}|)^{2/(2+\tau)}\over {h^{2}\eps^{2}}}\,
+\eps^{2/(2+\tau)}\Big), \nonumber\\
&=& \frac{1}{(n-i_p^\Delta+i_1^\Delta)h^{2p}\Delta} O\Big({\Delta^{1/(2+\tau)}\over
{h^{2}\eps^{2}}}\,
+\eps^{2/(2+\tau)}\Big) \nonumber\\
&=& \frac{1}{(n-i_p^\Delta+i_1^\Delta)h^{2p}\Delta} O\Big({\Delta^{(1-q)/2}\over
{h^{2}\eps^{2}}}\, +\eps^{1-q }\Big).
\end{eqnarray*}
Hence the last term in~(\ref{eq:varcov}) now gives the third order
bound (\ref{multvbound2a}). \hfill $\square$


\section*{References}

\small

\begin{verse}

Abramowitz, M.\ and Stegun, I\ (1964) , {\em Handbook of
Mathematical Functions, ninth edition}, Dover, New York.






Comte, F. \ (2004), Kernel deconvolution of stochastic volatility models, {\em J. Time Ser. Anal.} {\bf
25}, 563--582.

Comte, F.,  Dedecker, J. and Taupin, M.L.\ (2008), Adaptive density estimation for general ARCH models,
 {\em Econometric Theory} {\bf
24}, 1628--1662.

Comte, F. and Genon-Catalot, V.\ (2006), Penalized projection
estimator for volatility density, {\em Scand. J. Statist.} {\bf
33(4)}, 875--893.

Deo, C.M.\ (1973), A note on empirical processes for strong mixing
processes, {\em Ann. Probab.} {\bf 1}, 870--875.

Doukhan, P.\ (1994), {\em Mixing, Properties and Examples},
Springer-Verlag.

Fan, J.\ (1991), On the optimal rates of convergence for
nonparametric deconvolution problems, {\em Ann. Statist.} {\bf
19}, 1257--1272.


Franke, J., H\"ardle, W. and Kreiss, J.P. (2003), Nonparametric
estimation in a stochastic volatility model, In: Recent Advances
and Trends in Nonparametric Statistics, M.G. Akritas and D.N.
Politis Eds, Elsevier.

Genon-Catalot, V., Jeantheau, T.\ and Lar\'edo, C.\ (1998), Limit
theorems for discretely observed stochastic volatility models,
{\em Bernoulli} {\bf 4}, 283--303.

Genon-Catalot, V., Jeantheau, T.\ and Lar\'edo, C.\ (1999),
Parameter estimation
 for discretely observed stochastic volatility models,
{\em Bernoulli} {\bf 5}, 855-872.

Genon-Catalot, V., Jeantheau, T.\ and Lar\'edo, C.\ (2000),
Stochastic volatility models as hidden Markov models and
statistical applications, {\em Bernoulli} {\bf 6}, 1051--1079.

Gihman, I.I. and  Skorohod A.V.\ (1972), {\em Stochastic
Differential Equations}, Springer.

H\"ardle, W. (1990), {\em Smoothing Techniques}, Springer Verlag,
New York.

Heston, S.L.\  (1993), A closed-form solution for options with
stochastic volatility with applications to Bond and Currency
options, {\em The Review of Financial Studies} {\bf 6} (2),
327--343.

Hewitt, E.\ and Stromberg K.\ (1965), {\em Real and Abstract
Analysis}, Springer Verlag, New York.

Ibragimov, I.A., and Linnik, Yu.V. (1971), {\em Independent and
stationary sequences of random variables}, Wolters-Noordhoff.

Karatzas, I.\  and  Shreve, S.E. (1991), {\em Brownian Motion and
Stochastic Calculus}, Springer Verlag, New York.

Masry, E.\ (1983), Probability density estimation from sampled
data, {\em IEEE Trans. Inform. Theory}  {\bf 29}, 696--709.



Masry, E.\ (1993), Strong consistency and rates for
deconvolution of multivariate densities of stationary processes,
{\em Stoc. Proc. and Appl.} {\bf 475}, 53--74.

Nualart, D.\ (1995), {\em The Malliavin calculus and related
topics}, Springer Verlag, New York.


Scott, D.W.(1992), {\em Multivariate density estimation. Theory,
practice, and visualization}, Wiley, New York.

Skorokhod, A.V.\ (1989), {\em Asymptotic Methods in the Theory of
Stochastic Differential Equations}, AMS.


Van Es, B., Spreij, P. and Van Zanten H.\ (2003), Nonparametric
volatility density estimation, {\em Bernoulli} {\bf 9}, 451--465.

Van Es, B., Spreij, P. and Van Zanten H.\ (2005), Nonparametric
volatility density estimation for discrete time models, {\em J.
Nonparametr. Stat.} {\bf 17}, 237--251.

Van Zanten, H. and Zareba, P.\ (2008),  A note on wavelet density
deconvolution for weakly dependent data, {\em  Stat. Inference
Stoch. Process. } {\bf 11}, 207--219.

Wand, M.P.\ (1998), Finite sample performance of deconvolving
kernel density estimators, {\em Statist.\ Probab.\ Lett.} {\bf
37}, 131--139.

Wand, M.P.\ and Jones, M.C.\ (1995), Kernel Smoothing, Chapman and Hall,
London.

Wiggins, J.\ B.\  (1987), Option valuation under stochastic
volatility, {\em Journal of Financial Economics} {\bf 19},
351--372.

Van Zanten, H. and Zareba, P.\ (2008),  A note on wavelet density
deconvolution for weakly dependent data, {\em  Stat. Inference
Stoch. Process. } {\bf 11}, 207--219.

\end{verse}

\end{document}